\documentclass[12pt]{amsart}
\usepackage{amssymb}
\usepackage{amsfonts}
\usepackage{amssymb,latexsym}
\usepackage{enumerate}
\usepackage{mathrsfs}
\usepackage{tikz}
\usepackage{verbatim}
\makeatletter
\@namedef{subjclassname@2010}{%
  \textup{2010} Mathematics Subject Classification}
\makeatother

  [2002/01/19 v2.2g %
    AMS font definitions%
  ]
\DeclareFontFamily{U}{euf}{}
\DeclareFontShape{U}{euf}{m}{n}{%
  <5><6><7><8><9>gen*eufm%
  <10><10.95><12><14.4><17.28><20.74><24.88>eufm10%
  }{}
\DeclareFontShape{U}{euf}{b}{n}{%
  <5><6><7><8><9>gen*eufb%
  <10><10.95><12><14.4><17.28><20.74><24.88>eufb10%
  }{}

  [2002/01/19 v2.2g %
    AMS font definitions%
  ]
\DeclareFontFamily{U}{msb}{}
\DeclareFontShape{U}{msb}{m}{n}{%
  <5><6><7><8><9>gen*msbm%
  <10><10.95><12><14.4><17.28><20.74><24.88>msbm10%
  }{}

  [2002/01/19 v2.2g %
    AMS font definitions%
  ]
\DeclareFontFamily{U}{msa}{}
\DeclareFontShape{U}{msa}{m}{n}{%
  <5><6><7><8><9>gen*msam%
  <10><10.95><12><14.4><17.28><20.74><24.88>msam10%
  }{}

\newtheorem{theorem}{Theorem}[section]

\newtheorem{proposition}[theorem]{Proposition}

\newtheorem{corollary}[theorem]{Corollary}

\theoremstyle{definition}
\newtheorem{remark}[theorem]{Remark}

\numberwithin{equation}{section} 

\begin{document}

\title[]
{Integral representation of Weil's elliptic functions}

\author{Su Hu}
\address{Department of Mathematics, South China University of Technology, Guangzhou, Guangdong 510640, China}
\email{hus04@mails.tsinghua.edu.cn}

\author{Min-Soo Kim}
\address{Center for General Education, Kyungnam University,
7(Woryeong-dong) kyungnamdaehak-ro, Masanhappo-gu, Changwon-si, Gyeongsangnam-do 631-701, Republic of Korea
}
\email{mskim@kyungnam.ac.kr}

%\thanks{*Corresponding author}

\begin{abstract}In this paper, we give a two dimensional analogue of the Euler-MacLaurin summation formula. By using this formula, we obtain an integral representation of Weil's elliptic functions which was introduced in the book ``Elliptic functions according to Eisenstein and Kronecker".
\end{abstract}

\subjclass[2000]{11F03,11M35,33E05}
\keywords{Integral representation, Elliptic function}

%\thanks{Received May 24, 2009}

\maketitle

% \linenumbers


%\begin{proof}[Proof of Lemma 4]
%The result follows in the same way as in the proof of Proposition 1 in \cite{AD2}.
%\end{proof}

\section{Introduction}
\label{Intro}
The Hurwitz-Lerch zeta function $\Phi(z,s,a)$ is defined by
\begin{equation}\label{l-zeta}
\begin{aligned}
\Phi(z,s,a)
=\sum_{n=0}^\infty\frac{z^{n}}{(a+n)^{s}}
\end{aligned}
\end{equation}
for
$$a\in\mathbb C\setminus\mathbb Z_0^-;\quad s\in\mathbb C\text{ when }|z|<1; \quad \Re(s)>k\;(k\in\mathbb N) \text{ when } |z|=1.$$
Here $\mathbb Z_0^-=\{0,-1,-2,\ldots\},$ $\mathbb N$ is the set of positive integers, and $\mathbb C$ is the set of complex numbers
(see \cite{SC}).

Letting $z=1$ in (\ref{l-zeta}), we have the Hurwitz zeta function
\begin{equation}~\label{Hurwitz}\zeta(s,a)=\sum_{n=0}^{\infty}\frac{1}{(a+n)^{s}}.\end{equation}
Letting $a=1$ in (\ref{Hurwitz}), we obtain the Riemann zeta function \begin{equation}~\label{Riemann}\zeta_{\mathbb{Q}}(s)=\sum_{n=1}^{\infty}\frac{1}{n^{s}}.\end{equation}

Let $P_1(x)=B_1(x-[x])=x-[x]-1/2$ be the first periodized Bernoulli polynomials, and $\{x\}=x-[x]$ be the fractional part of $x.$ Recently, by using  Euler-MacLaurin summation formula \begin{equation}\label{du-sum-0}
\sum_{\alpha<n\leq\beta}\Phi(n)=\int_{\alpha}^\beta\Phi(x)dx +\int_{\alpha}^\beta\Phi'(x) P_1(x)dx \\
+P_1(\alpha)\Phi(\alpha)-P_1(\beta)\Phi (\beta),
\end{equation}
Coffey~\cite[p. 81]{JNT} gave an integral representation of the Hurwitz-Lerch zeta function $\Phi(z,s,a)$ (\ref{l-zeta}) as  follows.
\begin{proposition}[Coffey~\cite{JNT}]\label{Coffey}
\begin{equation}~\label{Cof}\begin{aligned}\Phi(z,s,a)&=\frac{1}{a^{s}}+\frac{z}{2(a+1)^{s}}\\&\quad+\int_{1}^{\infty}\frac{z^{x}}{(x+a)^{s}}dx+\int_{1}^{\infty}\left[\frac{z^{x}\ln z}{(x+a)^{s}}-\frac{sz^{x}}{(x+a)^{s+1}}\right]P_{1}(x)dx.\end{aligned}\end{equation}\end{proposition}

Suppose $W$ is a lattice in the complex plane, $\omega_{1}$ and $\omega_{2}$ are two generators of $W$, so that $W$ consists of the points $w=n\omega_{1}+m\omega_{2}$, where $n$ and $m$ are integers. In his historical book ``Elliptic functions according to Eisenstein and Kronecker"~\cite{Weil}, generalizing the Hurwitz zeta functions~(\ref{Hurwitz}) above, A. Weil~\cite[p. 14]{Weil} introduced the following elliptic function:
\begin{equation}~\label{weil}
E_{k}(a,W)=\sum_{w\in W}\frac{1}{(a+w)^{k}}=\sum_{(n,m)\in\mathbb{Z}^{2}}\frac{1}{(a+n\omega_{1}+m\omega_{2})^{k}}
\end{equation}
for $a\not\in W$, which is also a generalization of the homogeneous Eisenstein series defined by
\begin{equation}~\label{Eisenstein}
G_{k}(W)=\sum_{0\neq w\in W}\frac{1}{w^{k}}=\sum_{(n,m)\in\mathbb{Z}^2\backslash(0,0)}\frac{1}{(n\omega_{1}+m\omega_{2})^{k}}.
\end{equation}

As pointed out by Weil~\cite[p. 14]{Weil}, the series ~(\ref{weil}) are absolutely convergent for $k\geq 3$. If $k=1$ and $k=2$, then Eisenstein makes use of a summatory  process which we shall call Eisenstein summation., that is,
\begin{equation}~\label{ES}
\sum_{(n,m)\in\mathbb{Z}^{2}}\frac{1}{(a+n\omega_{1}+m\omega_{2})^{k}}=\lim_{M\to\infty}\sum_{m=-M}^{m=M}\left(\lim_{N\to\infty}\sum_{n=-N}^{N}\frac{1}{(a+n\omega_{1}+m\omega_{2})^{k}}\right).
\end{equation}

In order to give an analogue of the above Coffey's result for Weil's elliptic functions (\ref{weil}), in this paper, we shall generalize Euler-MacLaurin summation formula
(\ref{du-sum-0}) to a two dimensional case.

\begin{proposition}[Two dimensional  summation formula]\label{EM2}
Suppose $\Phi(x,y)=f(a+x\omega_{1}+y\omega_{2})\in C^{2}([\alpha_{1},\beta_{1}]\times[\alpha_{2},\beta_{2}])$, we have
\begin{equation}\label{main}
\sum_{\alpha_2<m\leq \beta_2}\sum_{\alpha_1<n\leq \beta_1} \Phi(n,m)
=I_1+I_2+I_3+I_4,
\end{equation}
where
\begin{equation}\label{main-eq1}
\begin{aligned}
I_1&=\iint_{[\alpha_1, \beta_1]\times[\alpha_{2},\beta_{2}]}\biggl[
\Phi(x,y) +\frac{\partial \Phi(x,y)}{\partial x} P_1(x) \\
&\quad +\frac{\partial \Phi(x,y)}{\partial y} P_1(y)
 +\frac{\partial^{2} \Phi(x,y)}{\partial x\partial y} P_1(x)P_1(y)
\biggl] dxdy,
\end{aligned}
\end{equation}

\begin{equation}\label{main-eq2}
\begin{aligned}
I_2&=\int_{\alpha_2}^{\beta_2}\biggl[
\Phi(\alpha_1,y)P_1(\alpha_1) - \Phi(\beta_1,y)P_1(\beta_1) \\
&\quad + \frac{\partial\Phi(\alpha_1,y)}{\partial y} P_1(y)P_1(\alpha_1)
-\frac{\partial \Phi(\beta_1,y)}{\partial y} P_1(y)P_1(\beta_1)
\biggl] dy ,
\end{aligned}
\end{equation}

\begin{equation}\label{main-eq3}
\begin{aligned}
I_3&=\int_{\alpha_1}^{\beta_1}\biggl[
\Phi(x,\alpha_2)P_1(\alpha_2)  - \Phi(x,\beta_2)P_1(\beta_2) \\
&\quad + \frac{\partial \Phi(x,\alpha_2)}{\partial x} P_1(x)P_1(\alpha_2)
-\frac{\partial \Phi(x,\beta_{2})}{\partial x} P_1(x)P_1(\beta_2)
\biggl] dx ,
\end{aligned}
\end{equation}
and
\begin{equation}\label{main-eq4}
\begin{aligned}
I_4&=P_1(\alpha_2)P_1(\alpha_1)\Phi(\alpha_1,\alpha_{2}) -P_1(\alpha_2)P_1(\beta_1)\Phi(\beta_1,\alpha_{2})  \\
&\quad-P_1(\beta_2)P_1(\alpha_1)\Phi(\alpha_1,\beta_{2}) +P_1(\beta_2)P_1(\beta_1)\Phi(\beta_1,\beta_{2}).
\end{aligned}
\end{equation}

\end{proposition}

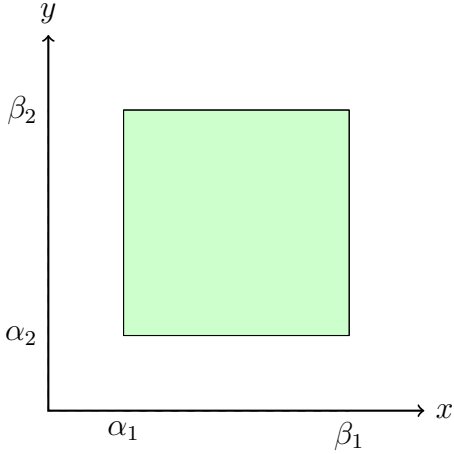
\begin{tikzpicture}\label{graph}
    \filldraw[fill=green!20](1,1) rectangle (4,4);
      \draw [<->,thick] (0,5) node (yaxis) [above] {$y$}
        |- (5,0) node (xaxis) [right] {$x$};
        \draw (1,1) coordinate (a) -- (1,4) coordinate (b);
     \draw (1,1) coordinate (a) -- (4,1) coordinate (c);;
     \draw (4,1) coordinate (c) -- (4,4) coordinate (d);
     \draw (1,4) coordinate (b) -- (4,4) coordinate (d);;
     \draw[dashed]   (xaxis -| a) node[below] {$\alpha_{1}$}
        -|(xaxis -| c) node[below] {$\beta_{1}$}
        -|(yaxis |- a) node[left]  {$\alpha_{2}$}
        -|(yaxis |- b) node[left] {$\beta_{2}$};;

\end{tikzpicture}

\begin{remark}From the above graph, we see that $I_1,$ the first term  in the right side of the above formula,
is an integral inside the square  $[\alpha_{1},\beta_{1}]\times[\alpha_{2},\beta_{2}]$. The second term $I_2$
and the third term $I_3$ are the integral along the boundary lines $[\alpha_{2},\beta_{2}]$ and $[\alpha_{1},\beta_{1}]$, respectively.
The fourth term
$I_4$ is a sum at the four corner points $(\alpha_1,\alpha_2), (\beta_1,\alpha_{2}), (\alpha_1,\beta_{2}), (\beta_1,\beta_{2})$.

Thus $I_1$ is corresponding to $\int_{\alpha}^\beta\Phi(x)dx +\int_{\alpha}^\beta\Phi'(x) P_1(x)dx$ (an integral inside the interval $[\alpha,\beta]$) in
the Euler-MacLaurin summation formula (see (\ref{du-sum-0}) above), while the term $I_2+I_3+I_4$ is corresponding to
the term $P_1(\alpha)\Phi(\alpha)-P_1(\beta)\Phi (\beta)$ (a sum at the boundary points of the interval $[\alpha,\beta]$) in
the Euler-MacLaurin summation formula (\ref{du-sum-0}).
\end{remark}

The following is an  integral representation of Weil's elliptic function (\ref{weil}) which is an analogue of Proposition~\ref{Coffey} on the integral representation of the Hurwitz-Lerch zeta function.

\begin{corollary}~\label{cor}
Suppose $(x_{0},y_{0})\in\mathbb R^2$ is the point such that $x_{0}w_1+y_{0}w_2=-a$, for arbitrary small $\epsilon>0$, we have
\begin{equation}~\label{weil2}
E_{k}(a,W)=\sum_{(n,m)\in\mathbb{Z}^{2}}\frac{1}{(a+n\omega_{1}+m\omega_{2})^{k}}
=J_1 +J_2+ J_3,
\end{equation}
where
\begin{equation}~\label{weil2-eq1}
\begin{aligned}
J_1&=\int_{-\infty}^{\infty}\biggl[
\frac1{(a+xw_1+(y_{0}+\epsilon)w_2)^k}P_1(y_{0}+\epsilon) \\
&\quad-\frac1{(a+xw_1+(y_{0}-\epsilon)w_2)^k}P_1(y_{0}-\epsilon) \\
&\quad + k\frac{w_1}{(a+xw_1+(y_{0}-\epsilon)w_2)^{k+1}}P_1(x)P_1(y_{0}-\epsilon)\\
&\quad - k\frac{w_1}{(a+xw_1+(y_{0}+\epsilon)w_2)^{k+1}}P_1(x)P_1(y_{0}+\epsilon)\biggl] dx,
\end{aligned}
\end{equation}

\begin{equation}~\label{weil2-eq2}
\begin{aligned}
J_2&=\iint_{[-\infty,\infty]\times[y_{0}+\epsilon,\infty]}\biggl[
\frac1{(a+xw_1+yw_2)^k} \\
&\quad-k\frac{w_1}{(a+xw_1+yw_2)^{k+1}}P_{1}(x) \\
&\quad-k\frac{w_2}{(a+xw_1+yw_2)^{k+1}}P_{1}(y)\\
&\quad+k(k+1)\frac{w_1w_2}{(a+xw_1+yw_2)^{k+2}}P_{1}(x)P_{1}(y)
\biggl]dxdy,
\end{aligned}
\end{equation}
and
\begin{equation}~\label{weil2-eq3}
\begin{aligned}
J_3&=\iint_{[-\infty,\infty]\times[-\infty,y_{0}-\epsilon]}\biggl[
\frac1{(a+xw_1+yw_2)^k} \\
&\quad-k\frac{w_1}{(a+xw_1+yw_2)^{k+1}}P_{1}(x) \\
&\quad-k\frac{w_2}{(a+xw_1+yw_2)^{k+1}}P_{1}(y)\\
&\quad+k(k+1)\frac{w_1w_2}{(a+xw_1+yw_2)^{k+2}}P_{1}(x)P_{1}(y)
\biggl]dxdy.
\end{aligned}
\end{equation}

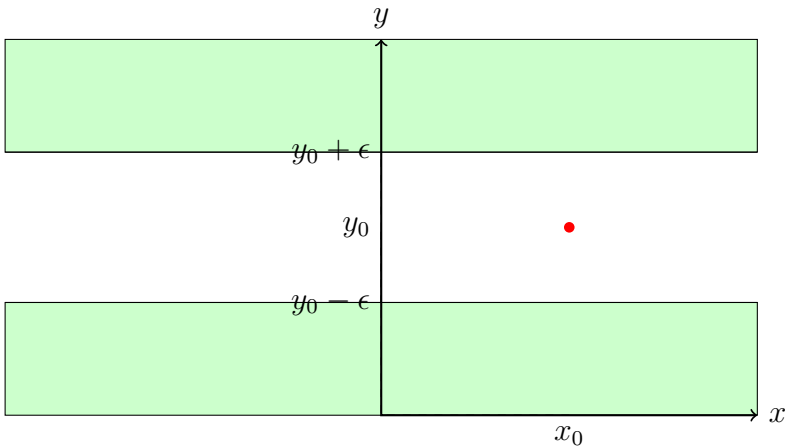
\begin{tikzpicture}\label{graph}

        \filldraw[fill=green!20](-5,3.5) rectangle (5,5);
         \filldraw[fill=green!20](-5,0) rectangle (5,1.5);
         \draw [<->,thick] (0,5) node (yaxis) [above] {$y$}
       |- (5,0) node (xaxis) [right] {$x$};
          \draw (-5,1.5) coordinate (a) -- (5,1.5) coordinate (c);;

     \draw (-5,3.5) coordinate (b) -- (5,3.5) coordinate (d);;
          \draw[dashed]     (xaxis -| (2.5,2.5) node[below] {$x_{0}$}
          -|(yaxis |- (2.5,2.5) node[left] {$y_{0}$}
        -|(yaxis |- a) node[left] {$y_{0}-\epsilon$}
        -|(yaxis |- b) node[left] {$y_{0}+\epsilon$};;
        \fill[red] (2.5,2.5) circle (2pt);
\end{tikzpicture}

\end{corollary}

\begin{remark} If $k=1$ and $k=2$, then the integral $\int_{-\infty}^\infty$ should be understood as $$\int_{-\infty}^\infty = \lim_{N\to\infty}\int_{-N}^{N}.$$
\end{remark}

In Sections 2 and 3, we shall prove Proposition ~\ref{EM2} and Corollary \ref{cor}, respectively.

\section{Proof of Proposition~\ref{EM2}}

The following Euler-MacLaurin summation formula will be used many times in our proof.

\begin{equation}\label{du-sum-01}
\sum_{\alpha<n\leq\beta}\Phi(n)=\int_{\alpha}^\beta\Phi(x)dx +\int_{\alpha}^\beta\Phi'(x) P_1(x)dx \\
+P_1(\alpha)\Phi(\alpha)-P_1(\beta)\Phi (\beta).
\end{equation}

By taking $\Phi(y)=\sum_{\alpha_1<n\leq \beta_1}f(a+nw_1+yw_2)$ in (\ref{du-sum-01}), we have
\begin{equation}\label{du-sum-1}
\begin{aligned}
\sum_{\alpha_2<m\leq \beta_2}\left(\sum_{\alpha_1<n\leq \beta_1}f(a+nw_1+mw_2)\right)
&=\int_{\alpha_2}^{\beta_2}\left(\sum_{\alpha_1<n\leq \beta_1}f(a+nw_1+yw_2)\right) dy \\
&\quad+\int_{\alpha_2}^{\beta_2}\left(\sum_{\alpha_1<n\leq \beta_1}\frac{\partial f(a+nw_1+yw_2)}{\partial y}\right)P_1(y) dy \\
&\quad+P_1(\alpha_2)\sum_{\alpha_1<n\leq \beta_1}f(a+nw_1+\alpha_2w_2) \\
&\quad-P_1(\beta_2)\sum_{\alpha_1<n\leq \beta_1}f(a+nw_1+\beta_2w_2).
\end{aligned}
\end{equation}

By taking $\Phi(x)=f(a+xw_1+yw_2)$ in (\ref{du-sum-01}), we have
\begin{equation}\label{du-sum-2}
\begin{aligned}
\sum_{\alpha_1<n\leq \beta_1}f(a+nw_1+mw_2)
&=\int_{\alpha_1}^{\beta_1}f(a+xw_1+yw_2) dx \\
&\quad+\int_{\alpha_1}^{\beta_1}\frac{\partial f(a+xw_1+yw_2)}{\partial x}P_1(x) dx \\
&\quad+P_1(\alpha_1)f(a+\alpha_1w_1+yw_2) \\
&\quad-P_1(\beta_1)f(a+\beta_1w_1+yw_2).
\end{aligned}
\end{equation}

By taking $\Phi(x)=\frac{\partial f(a+xw_1+yw_2)}{\partial y}$ in (\ref{du-sum-01}), we have
\begin{equation}\label{du-sum-3}
\begin{aligned}
\sum_{\alpha_1<n\leq \beta_1}\frac{\partial f(a+nw_1+yw_2)}{\partial y}
&=\int_{\alpha_1}^{\beta_1}\frac{\partial f(a+xw_1+yw_2)}{\partial y}dx  \\
&\quad+\int_{\alpha_1}^{\beta_1}\frac{\partial^{2} f(a+xw_1+yw_2)}{\partial x\partial y}P_1(x)dx \\
&\quad +P_1(\alpha_1)\frac{\partial f(a+\alpha_1w_1+yw_2)}{\partial y} \\
&\quad -P_1(\beta_1)\frac{\partial f(a+\beta_1w_1+yw_2)}{\partial y}.
\end{aligned}
\end{equation}

By taking $\Phi(x)=f(a+xw_1+\alpha_2w_2)$ in (\ref{du-sum-01}), we have
\begin{equation}\label{du-sum-4}
\begin{aligned}
\sum_{\alpha_1<n\leq \beta_1}f(a+nw_1+\alpha_2w_2)
&=\int_{\alpha_1}^{\beta_1}f(a+xw_1+\alpha_2w_2) dx \\
&\quad+\int_{\alpha_1}^{\beta_1}\frac{\partial f(a+xw_1+\alpha_2w_2)}{\partial x}P_1(x) dx \\
&\quad+P_1(\alpha_1)f(a+\alpha_1w_1+\alpha_2w_2) \\
&\quad-P_1(\beta_1)f(a+\beta_1w_1+\alpha_2w_2).
\end{aligned}
\end{equation}

By taking $\Phi(x)=f(a+xw_1+\beta_2w_2)$ in (\ref{du-sum-01}), we have
\begin{equation}\label{du-sum-5}
\begin{aligned}
\sum_{\alpha_1<n\leq \beta_1}f(a+nw_1+\beta_2w_2)
&=\int_{\alpha_1}^{\beta_1}f(a+xw_1+\beta_2w_2) dx \\
&\quad+\int_{\alpha_1}^{\beta_1}\frac{\partial f(a+xw_1+\beta_2w_2)}{\partial x}P_1(x) dx \\
&\quad+P_1(\alpha_1)f(a+\alpha_1w_1+\beta_2w_2) \\
&\quad-P_1(\beta_1)f(a+\beta_1w_1+\beta_2w_2).
\end{aligned}
\end{equation}

Substituting (\ref{du-sum-2}), (\ref{du-sum-3}), (\ref{du-sum-4}), (\ref{du-sum-5}) into (\ref{du-sum-1}), we have
$$
\begin{aligned}
\sum_{\alpha_2<m\leq \beta_2}\sum_{\alpha_1<n\leq \beta_1} f(a+nw_1+mw_2)
&= \sum_{\alpha_2<m\leq \beta_2}\left(\sum_{\alpha_1<n\leq \beta_1} f(a+nw_1+mw_2)\right) \\
&=\int_{\alpha_2}^{\beta_2}\biggl[
\int_{\alpha_1}^{\beta_1}f(a+xw_1+yw_2) dx \\
&\quad+\int_{\alpha_1}^{\beta_1}\frac{\partial f(a+xw_1+yw_2)}{\partial x}P_1(x) dx \\
&\quad+P_1(\alpha_1)f(a+\alpha_1w_1+yw_2) \\
&\quad-P_1(\beta_1)f(a+\beta_1w_1+yw_2)
\biggl]dy \\
&+\int_{\alpha_2}^{\beta_2}\biggl[
\int_{\alpha_1}^{\beta_1}\frac{\partial f(a+xw_1+yw_2)}{\partial y}dx  \\
&\quad+\int_{\alpha_1}^{\beta_1}\frac{\partial^{2}  f(a+xw_1+yw_2)}{\partial x\partial y}P_1(x)dx \\
&\quad +P_1(\alpha_1)\frac{\partial f(a+\alpha_1w_1+yw_2)}{\partial y} \\
&\quad-P_1(\beta_1)\frac{\partial f(a+\beta_1w_1+yw_2)}{\partial y}
\biggl]P_1(y)dy\\
&+ P_1(\alpha_2)\biggl[
\int_{\alpha_1}^{\beta_1}f(a+xw_1+\alpha_2w_2) dx \\
&\quad+\int_{\alpha_1}^{\beta_1}\frac{\partial f(a+xw_1+\alpha_2w_2)}{\partial x}P_1(x) dx \\
&\quad+P_1(\alpha_1)f(a+\alpha_1w_1+\alpha_2w_2) \\
&\quad-P_1(\beta_1)f(a+\beta_1w_1+\alpha_2w_2)
\biggl] \\
&- P_1(\beta_2)\biggl[
\int_{\alpha_1}^{\beta_1}f(a+xw_1+\beta_2w_2) dx \\
&\quad+\int_{\alpha_1}^{\beta_1}\frac{\partial f(a+xw_1+\beta_2w_2)}{\partial x}P_1(x) dx \\
&\quad+P_1(\alpha_1)f(a+\alpha_1w_1+\beta_2w_2) \\
&\quad-P_1(\beta_1)f(a+\beta_1w_1+\beta_2w_2)
\biggl]
\end{aligned}$$
This equals to
$$\begin{aligned}
&\int_{\alpha_2}^{\beta_2}\int_{\alpha_1}^{\beta_1}f(a+xw_1+yw_2) dxdy\\
&\quad
+\int_{\alpha_2}^{\beta_2}\int_{\alpha_1}^{\beta_1} \frac{\partial f(a+xw_1+yw_2)}{\partial x} P_1(x) dxdy \\
&\quad+ P_1(\alpha_1)\int_{\alpha_2}^{\beta_2}f(a+\alpha_1w_1+yw_2)dy\\
&\quad
- P_1(\beta_1)\int_{\alpha_2}^{\beta_2}f(a+\beta_1w_1+yw_2)dy   \\
&\quad+\int_{\alpha_2}^{\beta_2}\int_{\alpha_1}^{\beta_1} \frac{\partial f(a+xw_1+yw_2)}{\partial y} P_1(y) dxdy \\
&\quad +\int_{\alpha_2}^{\beta_2}\int_{\alpha_1}^{\beta_1} \frac{\partial^{2} f(a+xw_1+yw_2)}{\partial x\partial y} P_1(x)P_1(y) dxdy \\
&\quad+P_1(\alpha_1)\int_{\alpha_2}^{\beta_2}\frac{\partial f(a+\alpha_1w_1+yw_2)}{\partial y} P_1(y)dy \\
&\quad -P_1(\beta_1)\int_{\alpha_2}^{\beta_2}\frac{\partial f(a+\beta_1w_1+yw_2)}{\partial y} P_1(y)dy\\
&\quad+P_1(\alpha_2)\int_{\alpha_1}^{\beta_1}f(a+xw_1+\alpha_2w_2)dx\\
&\quad
+P_1(\alpha_2)\int_{\alpha_1}^{\beta_1}\frac{\partial f(a+xw_1+\alpha_2w_2)}{\partial x} P_1(x)dx \\
&\quad+P_1(\alpha_2)P_1(\alpha_1)f(a+\alpha_1w_1+\alpha_2w_2)\\
&\quad
-P_1(\alpha_2)P_1(\beta_1)f(a+\beta_1w_1+\alpha_2w_2) \\
&\quad-P_1(\beta_2)\int_{\alpha_1}^{\beta_1}f(a+xw_1+\beta_2w_2)dx\\
&\quad
-P_1(\beta_2)\int_{\alpha_1}^{\beta_1}\frac{\partial f(a+xw_1+\beta_2w_2)}{\partial x} P_1(x)dx \\
&\quad-P_1(\beta_2)P_1(\alpha_1)f(a+\alpha_1w_1+\beta_2w_2)\\
&\quad
+P_1(\beta_2)P_1(\beta_1)f(a+\beta_1w_1+\beta_2w_2).
\end{aligned}
$$

Finally, we have
\begin{equation}\label{fl}
\begin{aligned}
\sum_{\alpha_2<m\leq \beta_2}\sum_{\alpha_1<n\leq \beta_1} f(a+nw_1+mw_2)
&=\int_{\alpha_2}^{\beta_2}\int_{\alpha_1}^{\beta_1}\biggl[
f(a+xw_1+yw_2) \\
&\quad +\frac{\partial f(a+xw_1+yw_2)}{\partial x} P_1(x) \\
&\quad +\frac{\partial f(a+xw_1+yw_2)}{\partial y} P_1(y) \\
&\quad +\frac{\partial^{2} f(a+xw_1+yw_2)}{\partial x\partial y} P_1(x)P_1(y)
\biggl] dxdy \\
&+\int_{\alpha_2}^{\beta_2}\biggl[
f(a+\alpha_1w_1+yw_2)P_1(\alpha_1)  \\
&\quad- f(a+\beta_1w_1+yw_2)P_1(\beta_1) \\
&\quad + \frac{\partial f(a+\alpha_1w_1+yw_2)}{\partial y} P_1(y)P_1(\alpha_1) \\
&\quad-\frac{\partial f(a+\beta_1w_1+yw_2)}{\partial y} P_1(y)P_1(\beta_1)
\biggl] dy  \\
&+\int_{\alpha_1}^{\beta_1}\biggl[
f(a+xw_1+\alpha_2w_2)P_1(\alpha_2)  \\
&\quad- f(a+xw_1+\beta_2w_2)P_1(\beta_2) \\
&\quad + \frac{\partial f(a+xw_1+\alpha_2w_2)}{\partial x} P_1(x)P_1(\alpha_2) \\
&\quad-\frac{\partial f(a+xw_1+\beta_2w_2)}{\partial x} P_1(x)P_1(\beta_2)
\biggl] dx \\
&\quad+P_1(\alpha_2)P_1(\alpha_1)f(a+\alpha_1w_1+\alpha_2w_2) \\
&\quad-P_1(\alpha_2)P_1(\beta_1)f(a+\beta_1w_1+\alpha_2w_2) \\
&\quad-P_1(\beta_2)P_1(\alpha_1)f(a+\alpha_1w_1+\beta_2w_2) \\
&\quad+P_1(\beta_2)P_1(\beta_1)f(a+\beta_1w_1+\beta_2w_2).
\end{aligned}
\end{equation}
Notice that $\Phi(x,y)=f(a+x\omega_{1}+y\omega_{2}),$ we get our result.

\section{Proof of Corollary~\ref{cor}}
In this section, let
\begin{equation}\label{du-sum-7}
f(a,x,y)=\frac1{(a+xw_1+yw_2)^k}.
\end{equation}
We have
\begin{equation}\label{du-sum-8}
\frac{\partial f}{\partial x}(x,y)=-k\frac{w_1}{(a+xw_1+yw_2)^{k+1}},
\end{equation}

\begin{equation}\label{du-sum-9}
\frac{\partial^{2}  f}{\partial x\partial y}(x,y)=k(k+1)\frac{w_1w_2}{(a+xw_1+yw_2)^{k+2}},
\end{equation}

\begin{equation}\label{du-sum-10}
\frac{\partial f}{\partial y}(x,y)=-k\frac{w_2}{(a+xw_1+yw_2)^{k+1}},
\end{equation}

\begin{equation}\label{du-sum-11}
\frac{\partial f}{\partial x}(N,y)=-k\frac{w_1}{(a+Nw_1+yw_2)^{k+1}}
\end{equation}
and
\begin{equation}\label{du-sum-11}
\frac{\partial f}{\partial y}(x,N)=-k\frac{w_2}{(a+xw_1+Nw_2)^{k+1}}.
\end{equation}
We have
\begin{equation}\label{du-sum-12}
\begin{aligned}
&|f(a,N,y)|=\frac1{|a+Nw_1+yw_2|^k}\rightarrow 0,\\
&\left|\frac{\partial f}{\partial x}(N,y)\right|=k\frac{|w_1|}{|a+Nw_1+yw_2|^{k+1}}\rightarrow 0 \quad\text{as } N\to \infty,
\end{aligned}
\end{equation}uniformly for $y\in [\alpha_{2},\beta_{2}]$,
and
\begin{equation}\label{du-sum-131}
\begin{aligned}
&|f(a,x,N)|=\frac1{|a+xw_1+Nw_2|^k}\rightarrow 0,\\
&\left|\frac{\partial f}{\partial y}(x,N)\right|=k\frac{|w_2|}{|a+xw_1+Nw_2|^{k+1}}\rightarrow 0 \quad\text{as } N\to \infty,
\end{aligned}
\end{equation}uniformly for $x\in [\alpha_{1},\beta_{1}]$.

Thus, by taking $\alpha_{1}=-N$, $\beta_{1}=N$, $\alpha_{2}=y_{0}+\epsilon$, $\beta_{2}=N$ and $f(a,x,y)=\frac1{(a+xw_1+yw_2)^k}$ in (\ref{fl}), then letting $N\to\infty$, we have
\begin{equation}\label{du-sum-15}
\begin{aligned}
\sum_{y_{0}+\epsilon<m<\infty}\sum_{-\infty <n<\infty} \frac1{(a+nw_1+mw_2)^k}
&=\int_{y_{0}+\epsilon}^\infty\int_{-\infty}^\infty\biggl[
\frac1{(a+xw_1+yw_2)^k} \\
&\quad-k\frac{w_1}{(a+xw_1+yw_2)^{k+1}}P_{1}(x) \\
&\quad-k\frac{w_2}{(a+xw_1+yw_2)^{k+1}}P_{1}(y)\\
&\quad+k(k+1)\frac{w_1w_2}{(a+xw_1+yw_2)^{k+2}}P_{1}(x)P_{1}(y)
\biggl]dxdy\\
&\quad+\int_{-\infty}^{\infty}\biggl[
\frac1{(a+xw_1+(y_{0}+\epsilon)w_2)^k}P_1(y_{0}+\epsilon) \\
&\quad - k\frac{w_1}{(a+xw_1+(y_{0}+\epsilon)w_2)^{k+1}}P_1(x)P_1(y_{0}+\epsilon)\biggl] dx .
\end{aligned}
\end{equation}
Similarly, we have
\begin{equation}\label{du-sum-16}
\begin{aligned}
\sum_{-\infty <m\leq y_{0}-\epsilon }\sum_{-\infty <n<\infty} \frac1{(a+nw_1+mw_2)^k}
&=\int_{-\infty}^{y_{0}-\epsilon}\int_{-\infty}^\infty\biggl[
\frac1{(a+xw_1+yw_2)^k} \\
&\quad-k\frac{w_1}{(a+xw_1+yw_2)^{k+1}}P_{1}(x) \\
&\quad-k\frac{w_2}{(a+xw_1+yw_2)^{k+1}}P_{1}(y)\\
&\quad+k(k+1)\frac{w_1w_2}{(a+xw_1+yw_2)^{k+2}}P_{1}(x)P_{1}(y)
\biggl]dxdy\\
&\quad-\int_{-\infty}^{\infty}\biggl[
\frac1{(a+xw_1+(y_{0}-\epsilon)w_2)^k}P_1(y_{0}-\epsilon) \\
&\quad + k\frac{w_1}{(a+xw_1+(y_{0}-\epsilon)w_2)^{k+1}}P_1(x)P_1(y_{0}-\epsilon)\biggl] dx .
\end{aligned}
\end{equation}
Since $a\not\in W$, for arbitrary small $\epsilon>0$, we have
\begin{equation}\label{du-sum-17}
\begin{aligned}\sum_{(n,m)\in\mathbb{Z}^{2}}\frac{1}{(a+n\omega_{1}+m\omega_{2})^{k}}&=\sum_{y_{0}+\epsilon<m<\infty}\sum_{-\infty <n<\infty}\frac1{(a+nw_1+mw_2)^k}\\
&\quad+\sum_{-\infty <m\leq y_{0}-\epsilon }\sum_{-\infty <n<\infty} \frac1{(a+nw_1+mw_2)^k},\end{aligned}\end{equation}
from (\ref{du-sum-15}) and (\ref{du-sum-16}), we get our result.

\bibliography{central}

\end{document}